\documentclass[11pt]{article}
\usepackage{amssymb,amsfonts,amsmath,curves}
\usepackage{epsfig}
\usepackage{graphicx}
\usepackage{color}
\usepackage{anysize}
\usepackage{hyperref}
\usepackage{url}
\usepackage{enumitem}
\marginsize{2cm}{2cm}{1cm}{1cm}

\begin{document}

\makeatletter\@addtoreset{equation}{section}
\makeatother\def\theequation{\thesection.\arabic{equation}}

\baselineskip=14pt
\parskip=4pt

\newtheorem{teo}{Theorem}[section]
\newtheorem{defin}[teo]{Definition}
\newtheorem{con}[teo]{Conjecture}
\newtheorem{prop}[teo]{Proposition}
\newtheorem{lem}[teo]{Lemma}
\newtheorem{rmk}[teo]{Remark}
\newtheorem{cor}[teo]{Corollary}
\newtheorem{clm}[teo]{Claim}

\newcommand{\1}{{\text{\Large $\mathfrak 1$}}}

\newcommand{\be}{\begin{equation}}
\newcommand{\ee}{\end{equation}}
\newcommand{\beqn}{\begin{eqnarray}}
\newcommand{\beqnn}{\begin{eqnarray*}}
\newcommand{\eeqn}{\end{eqnarray}}
\newcommand{\eeqnn}{\end{eqnarray*}}
\newcommand{\bp}{\begin{prop}}
\newcommand{\ep}{\end{prop}}
\newcommand{\bt}{\begin{teo}}
\newcommand{\bcor}{\begin{cor}}
\newcommand{\ecor}{\end{cor}}
\newcommand{\bcon}{\begin{con}}
\newcommand{\econ}{\end{con}}
\newcommand{\bconj}{\begin{conj}}
\newcommand{\econj}{\end{conj}}
\newcommand{\et}{\end{teo}}
\newcommand{\brm}{\begin{rmk}}
\newcommand{\erm}{\end{rmk}}
\newcommand{\blem}{\begin{lem}}
\newcommand{\elem}{\end{lem}}
\newcommand{\ben}{\begin{enumerate}}
\newcommand{\een}{\end{enumerate}}
\newcommand{\bei}{\begin{itemize}}
\newcommand{\eei}{\end{itemize}}
\newcommand{\bdf}{\begin{defin}}
\newcommand{\edf}{\end{defin}}
\newcommand{\bpr}{\begin{proof}}
\newcommand{\epr}{\end{proof}}

\newenvironment{proof}{\noindent {\em Proof}.\,\,}{\hspace*{\fill}$\halmos$\medskip}

\newcommand{\halmos}{\rule{1ex}{1.4ex}}
\def \qed {{\hspace*{\fill}$\halmos$\medskip}}

\newcommand{\fr}{\frac}
\newcommand{\Z}{{\mathbb Z}}
\newcommand{\R}{{\mathbb R}}
\newcommand{\E}{{\mathbb E}}
\newcommand{\C}{{\mathbb C}}
\renewcommand{\P}{{\mathbb P}}
\newcommand{\N}{{\mathbb N}}
\newcommand{\Q}{{\mathbb Q}}
\newcommand{\U}{{\mathbb U}}
\newcommand{\var}{{\mathbb V}}
\renewcommand{\S}{{\cal S}}
\newcommand{\T}{{\cal T}}
\newcommand{\W}{{\cal W}}
\newcommand{\X}{{\cal X}}
\newcommand{\Y}{{\cal Y}}
\newcommand{\h}{{\cal H}}
\newcommand{\Bi}{{\cal B}}
\newcommand{\Di}{{\cal D}}
\newcommand{\Ei}{{\cal E}}
\newcommand{\Fi}{{\cal F}}
\newcommand{\Li}{{\cal L}}
\newcommand{\Mi}{{\cal M}}
\newcommand{\Ni}{{\cal N}}
\newcommand{\Ui}{{\cal U}}

\renewcommand{\a}{\alpha}
\renewcommand{\b}{\beta}
\newcommand{\g}{\gamma}
\newcommand{\G}{\Gamma}
\renewcommand{\L}{\Lambda}
\renewcommand{\l}{\lambda}
\renewcommand{\d}{\delta}
\newcommand{\D}{\Delta}
\newcommand{\e}{\epsilon}
\newcommand{\eps}{\epsilon}
\newcommand{\s}{\sigma}
\newcommand{\B}{{\cal B}}
\renewcommand{\o}{\omega}

\newcommand{\nn}{\nonumber}
\renewcommand{\=}{&=&}
\renewcommand{\>}{&>&}
\newcommand{\<}{&<&}
\renewcommand{\le}{\leq}
\newcommand{\+}{&+&}

\newcommand{\pa}{\partial}
\newcommand{\ffrac}[2]{{\textstyle\frac{{#1}}{{#2}}}}
\newcommand{\dif}[1]{\ffrac{\partial}{\partial{#1}}}
\newcommand{\diff}[1]{\ffrac{\partial^2}{{\partial{#1}}^2}}
\newcommand{\difif}[2]{\ffrac{\partial^2}{\partial{#1}\partial{#2}}}
\newcommand{\asto}[1]{\underset{{#1}\to\infty}{\longrightarrow}}
\newcommand{\Asto}[1]{\underset{{#1}\to\infty}{\Longrightarrow}}

\definecolor{Red}{rgb}{1,0,0}

\newcommand{\red}[1]{{\color{red}{#1}}}
\newcommand{\blue}[1]{{\color{blue}{#1}}}
\newcommand{\Ex}[1]{\mathbb{E}\!\left[#1\right]}
\newcommand{\pr}[1]{\mathbb{P}\!\left(#1\right)}
\newcommand{\vol}[1]{\mathrm{Vol}\!\left(#1\right)}

\title{Symmetric Rearrangements Around Infinity with Applications to L\'evy Processes}
\author{Alexander Drewitz$^{\,1}$ \and Perla Sousi$^{\,2}$ \and Rongfeng Sun$^{\,3}$}
\date{Jan 28, 2013}
\maketitle

\footnotetext[1]{Department of Mathematics, Columbia University, New York, USA. Email: drewitz@math.columbia.edu}

\footnotetext[2]{Statistical Laboratory, University of Cambridge, Cambridge, UK. Email: p.sousi@statslab.cam.ac.uk}

\footnotetext[3]{Department of Mathematics, National University of Singapore, Singapore. Email: matsr@nus.edu.sg}

\begin{abstract}
We prove a new rearrangement inequality for multiple integrals, which partly generalizes a result of Friedberg and Luttinger~\cite{FL76} and can be interpreted as involving symmetric rearrangements of domains around $\infty$. As applications, we prove two comparison results for general L\'evy processes and their symmetric rearrangements. The first application concerns the survival probability of a point particle in a Poisson field of moving traps following independent L\'evy motions. We show that the survival probability can only increase if the point particle does not move, and the traps and the L\'evy motions are symmetrically rearranged. This essentially generalizes an isoperimetric inequality of Peres and Sousi~\cite{PS11} for the Wiener sausage. In the second application, we show that the $q$-capacity of a Borel measurable set for a L\'evy process can only decrease if the set and the L\'evy process are symmetrically rearranged. This result generalizes an inequality obtained by Watanabe~\cite{W83} for symmetric L\'evy processes.


\medskip

\noindent
\emph{Keywor{\rm d}s:} capacity, isoperimetric inequality, L\'evy process, L\'evy sausage, Pascal principle, rearrangement inequality,  trapping {\rm d}ynamics.
\end{abstract}

\section{Introduction}
\subsection{Rearrangement Inequality}
As motivation, let us start with the following random walk exit problem. Suppose that $(X_n)_{n\geq 0}$ is a discrete time random walk on $\R^d$ with transition probability kernel $p_{n}(x){\rm d}x$ from time $n-1$ to $n$. Let $(A_n)_{n\geq 0}$ be a sequence of Borel-measurable sets in $\R^d$ with finite volume, such that the walk is killed at time $i$ if $X_i\notin A_i$. If $X_0$ is uniformly distributed on $A_0$, then
\be\label{survprob}
\P(X_i \in A_i \ \forall\ 0\leq i\leq n) = \frac{1}{|A_0|} \idotsint \prod_{i=0}^n 1_{A_i}(x_i) \prod_{i=1}^n p_i(x_i-x_{i-1}) \prod_{i=0}^n {\rm d}x_i,
\ee
where $|A_0|$ denotes the Lebesgue measure of $A_0$. By the classic Brascamp-Lieb-Luttinger rearrangement inequality (see \cite{BLL74} and~\cite[Theorem 3.8]{LL01}), the above probability is upper bounded by
\be
\frac{1}{|A^*_0|} \idotsint \prod_{i=0}^n 1_{A^*_i}(x_i) \prod_{i=1}^n p^*_i(x_i-x_{i-1}) \prod_{i=0}^n {\rm d}x_i,
\ee
where $A_i^*$ and $p_i^*$ denote respectively the symmetric decreasing rearrangements of $A_i$ and $p_i$, which are defined as follows.
\bdf\label{D:rearrange}\rm{
If $A\subset \R^d$ with $|A|<\infty$ (i.e., $A$ has finite volume), then its symmetric decreasing rearrangement $A^*$ is defined to be the open ball centered at the origin with $|A^*|=|A|$. If $|A|=\infty$, then we define $A^* := \R^d$. If $f: \R^d\to [0,\infty]$ is measurable, then its symmetric decreasing rearrangement $f^*$ is defined to
be
\[
f^*(x) := \int_0^{\infty}1_{F_t^*}(x) \, {\rm d}t, \qquad x\in\R^d,
\]
where $F_t := \{y: f(y)> t\}$, $t\geq 0$, are the level sets of $f$ (note that $f(x) = \int_{0}^{\infty} 1_{F_t}(x) \,{\rm d}t$). In particular, $f^*(x)=g(|x|)$ for a $g:[0,\infty)\to [0,\infty]$ which is nonincreasing and right-continuous. }
\edf
In other wor{\rm d}s, in (\ref{survprob}), the probability that the walk $X$ survives up to time $n$ can only increase if its transition kernels, as well as the domains, are all replaced by their symmetric decreasing rearrangements. There is a sizable literature on rearrangement inequalities and their relation to isoperimetric problems, see e.g.~\cite[Chapter 3]{LL01}. Combined with probabilistic representations, rearrangement inequalities can be used to obtain the celebrated Rayleigh-Faber-Krahn inequality on the first eigenvalue of the Dirichlet Laplacian, and comparison inequalities for heat kernels and Green functions (see e.g.~\cite{BS01, BM-H10} and the references therein).

We are interested in the analogue of the survival probability in (\ref{survprob}), where we replace the domains $A_i$ by their complements $A_i^c$. Since
$|A_i^c|=\infty$, the multiple integral in (\ref{survprob}) is in general infinite if we replace $A_i$ by $A_i^c$. However, it is sensible to consider
instead
\be\label{exitmass}
W_n((A_i)_{i\geq 0}, (p_i)_{i\geq 1}):=\idotsint \Big(1-\prod_{i=0}^n 1_{A^c_i}(x_i)\Big) \prod_{i=1}^n p_i(x_i-x_{i-1}) \prod_{i=0}^n {\rm d}x_i.
\ee
If we interpret $\prod_{i=1}^n p_i(x_i-x_{i-1}) \prod_{i=0}^n {\rm d}x_i$ as an infinite measure on the space of random walk trajectories, then $W_n$ can be interpreted as the total measure of the trajectories of $X$ that are killed by the hard traps $(A_i)_{i\geq 0}$ by time $n$, where the initial measure of $X_0$ is the Lebesgue measure on $\R^d$ instead of a probability measure. The quantity $W_n$ is also equal to the expected volume of the ``sausage" based on the sets $(-A_i)_{i\leq n}$ around the walk $X$ with $X_0=0$ and transition density $p_i(x)$, that is,
\be
W_n((A_i)_{i\geq 0}, (p_i)_{i\geq 1}) = \Ex{\vol{\cup_{i=0}^{n} (X_i - A_i)}}.
\ee
The rearrangement inequality we will prove amounts to the statement that
\be\label{ri0}
W_n((A_i)_{i\geq 0}, (p_i)_{i\geq 1}) \geq W_n((A^*_i)_{i\geq 0}, (p^*_i)_{i\geq 1}).
\ee
Although (\ref{ri0}) is still formulated in terms of symmetric decreasing rearrangements of $A_i$ and $p_i$, with the origin being the center of rearrangements, it does not follow directly from classic rearrangement inequalities because terms with alternating signs appear when we expand $\prod_{i=0}^n (1-1_{A_i}(x_i))$. In both
(\ref{survprob}) and (\ref{exitmass}), the goal is to maximize the probability that the walk stays within the domains. The only difference is the replacement of
the domains $(A_i)_{i\geq 1}$ in (\ref{survprob}) by their complements in (\ref{exitmass}). In light of the close analogy between the two problems, it is instructive to think of (\ref{ri0}) as a rearrangement inequality where the infinite domains $A_i^c$ are symmetrically rearranged around $\infty$. This point of view will guide our proof.

We now formulate our rearrangement inequality for multiple integrals, which is a more general version of (\ref{ri0}). We will assume that: The initial measure for $X_0$
is $\phi(x){\rm d}x$ for some $\phi: \R^d\to [0,\infty)$; each hard trap $A_i$ is replaced by a trap function $V_i: \R^d\to [0,1]$, so that upon jumping to $x_i$ at time $i$, the walk is killed with probability $V_i(x_i)$ instead of $1_{A_i}(x_i)$; each kernel $p_i: \R^d\to [0,\infty)$ is no longer assumed to be a probability density kernel.

\bt\label{T:ri}
Let $\phi: \R^d \to [0,\infty)$ and let $\sigma:=\sup\{ t\geq 0: |\{x:\phi(x)<t\}|<\infty\}$. Define the symmetric increasing rearrangement of $\phi$ by
$\phi_*:=\sigma-(\sigma-\phi\wedge\sigma)^*$. For $i\geq 0$ and $j\geq 1$, let $V_i: \R^d\to [0,1]$ and $p_j: \R^d\to [0,\infty)$, and let $V_i^*$ and $p_j^*$ denote their symmetric decreasing rearrangements. Denote $V_\cdot:=(V_i)_{i\geq 0}$, $p_\cdot:=(p_j)_{j\geq 1}$, $V_\cdot^*:=(V_i^*)_{i\geq 0}$, and $p_\cdot^*:=(p_j^*)_{j\geq 1}$. Then for all $n\geq 0$,
\be\label{ri}
\begin{aligned}
&\ W_n(\phi, V_\cdot, p_\cdot) \!\!\! &:=& \idotsint  \phi(x_0) \Big(1-\prod_{i=0}^n (1-V_i(x_i))\Big) \prod_{i=1}^n p_i(x_i-x_{i-1})  \prod_{i=0}^n {\rm d}x_i \\
\geq\  &\ W_n(\phi_*, V^*_\cdot, p^*_\cdot) \!\!\! &:=&  \idotsint   \phi_*(x_0)  \Big(1-\prod_{i=0}^n(1-V^*_i(x_i))\Big) \prod_{i=1}^n p^*_i(x_i-x_{i-1}) \prod_{i=0}^n {\rm d}x_i.
\end{aligned}
\ee
\et
\brm\label{R:FL}{\rm We will discuss two extensions of (\ref{ri}) in Remark~\ref{R:riext}.
Theorem~\ref{T:ri} partly generalizes an inequality of Friedberg and Luttinger~\cite[Corollary~2]{FL76}, which is the special case of (\ref{ri}) in dimension $d=1$, with $\phi\equiv 1$ and $p_i=p_i^*$ for all $i\geq 1$. They however allow an additional convolution kernel $p_{n+1}(x_0-x_n)$ with $p_{n+1}=p_{n+1}^*$, which is set to $1$ in our case. As pointed out at the end of~\cite{FL76}, if we include the additional kernel $p_{n+1}(x_0-x_n)$, then the analogue of (\ref{ri}) is false in general. Indeed, if we let $n=1$, $\phi\equiv 1$, $V_0\equiv V_1\equiv 1$, and choose $p_1$ and $p_2$ such that $(p_1*p_2)(0)=0$ and $(p_1^* * p_2^*)(0)>0$, then the analogue of (\ref{ri}) reads as
$$
\iint p_1(x_1-x_0) p_2(x_0-x_1) \, {\rm d}x_1 {\rm d} x_0 = \int (p_1*p_2)(0) \, {\rm d} x_0 =0 \geq \int (p_1^* *p_2^*)(0) \, {\rm d} x_0 =\infty,
$$
which is clearly false.}
\erm

Although stated only for dimension $1$ in~\cite{FL76}, Friedberg and Luttinger's inequality extends to higher dimensions by standard symmetrization techniques developed in~\cite{BLL74}, as noted in~\cite{M-H06}. Recently, Peres and Sousi~\cite[Prop.~1.6]{PS11} gave a different proof of this fact. More precisely, they proved (\ref{ri}) where $(V_i)_{i\geq 0}$ were taken to be indicator functions of open sets, and $(p_i)_{\geq 1}$ were taken to be the densities of uniform distributions on centered open balls. The interpretation of symmetric rearrangements around $\infty$ arises naturally in their proof. They appealed to an analogue rearrangement inequality on the sphere by Burchard and Schmuckenschl\"ager~\cite[Theorem 2]{BS01}, which they applied by performing symmetric decreasing rearrangements of domains around the south pole of the sphere. As the radius of the sphere ten{\rm d}s to infinity, the neighborhood around the north pole approximates $\R^d$, while the south pole converges to $\infty$.
When we consider the case which requires symmetric decreasing rearrangements of the convolution kernels $(p_i)_{i\geq 1}$ as in (\ref{ri}), there appear to be no existing analogous rearrangement inequalities on the sphere that we can appeal to. Instead, we develop a more direct and surprisingly simple approach to prove (\ref{ri}). Our proof contains two ingredients. The first is induction over the number of factors (or time steps) in the integrands in (\ref{ri}), which makes essential use of the Markovian structure of the problem. Similar induction approaches to rearrangement inequalities have been used before, see e.g.~\cite{B94}. The second and key ingredient is a proper notion of symmetric domination, which is motivated by the point of view of symmetric rearrangements around $\infty$. Readers who are interested in the proof of Theorem~\ref{T:ri} can jump directly to Section~\ref{S:ri}, where the simple proof is presented.
\bigskip

Our primary motivation for Theorem~\ref{T:ri} originates in the stu{\rm d}y of the survival probability of a point particle in a Poisson field of moving traps, each following an independent L\'evy motion, which gives rise to continuous time analogues of the total killed measure $W_n$ defined in (\ref{exitmass}). Our first application of Theorem~\ref{T:ri} is to show that the survival probability of the point particle can only increase if it stays put, while the L\'evy motions and the shape of the traps are symmetric decreasingly rearranged (see Theorem \ref{T:trap}). Previously, Peres and Sousi~\cite{PS11} proved such a comparison result when the traps follow independent Brownian motions, so that only the point particle motion and the shape of the traps require symmetric decreasing rearrangements. Our attempt to generalize their result to allow for symmetric decreasing rearrangements of general L\'evy motions was inspired by the work of Ba\~nuelos and M\'endez-Hern\'andez~\cite{BM-H10}, where a continuous time analogue of the exit problem in (\ref{survprob}) was considered. More specifically, they showed that the survival probability of a L\'evy motion in a time-independent trap potential on a finite volume open domain can only increase if the L\'evy motion and the domain are symmetric decreasingly rearranged, while the trap potential is symmetric increasingly rearranged\footnote{There was an error in the formulation of Theorem 1.4 in~\cite{BM-H10}, where the symmetric decreasing rearrangement $V^*$ of the potential $V$ on
the domain $D$ should be replaced by its symmetric increasing rearrangement $V_*$ on the domain $D^*$.}.

Like classical rearrangement inequalities, Theorem~\ref{T:ri} also has its potential-theoretic implications. As a second application of Theorem~\ref{T:ri}, we prove a comparison inequality for capacities of sets for L\'evy processes (Theorem \ref{T:cap}). More precisely, we show that if $A$ is any Borel-measurable subset of $\R^d$, then the $q$-capacity of $A$ for a L\'evy process $X$ ($q>0$ if $X$ is recurrent, and $q\geq 0$ if $X$ is transient) can only decrease if we replace $A$ and $X$ by their symmetric decreasing rearrangements. This generalizes a result of Watanabe~\cite{W83}, who proved such a comparison inequality
for symmetric L\'evy processes using Dirichlet forms. Special cases of Watanabe's result have been reproduced by Betsakos~\cite{B04} and M\'endez-Hern\'andez~\cite{M-H06}. An inequality of the type in Theorem~\ref{T:ri} was in fact conjectured in~\cite{BM-H10}, where its connection to $0$-capacities was also pointed out.

In the remainder of this introduction, we will formulate precisely our comparison inequalities for the trapping problem and for capacities. We will then end the
introduction with an outline of the rest of the paper.

\subsection{Trapping Problem}\label{S:IntroTrap}
The model of a point particle in a Poisson field of moving traps in $\R^d$ is defined as follows. The point particle follows a deterministic path in $\R^d$, given
by the function $f: [0,\infty)\to\R^d$. Let $\Xi_0$ be a Poisson point process on $\R^d$ with intensity measure $\phi(x){\rm d}x$ for some $\phi: \R^d\to [0,\infty)$. We label the points in $\Xi_0$ by $(z^n_0)_{n\in\N}$. The points in $\Xi_0$ move independently in time, each following the law of a L\'evy process $X:=(X_t)_{t\geq 0}$. Namely, we replace $\Xi_0$ at time $t>0$ by $\Xi_t:=\{z^n_t : n\in\N\}$, where $z^n_t = z^n_0+ X^n_t$, and $(X^n)_{n\in\N}$ are i.i.d.\ copies of $X$ with $X_0=0$. The points in $\Xi_t$ determine the location of traps at time $t$, and the actual shape of the traps at time $t$ is determined by a trap potential $U_t : \R^d\to [0,\infty]$. More precisely, the field of traps at time $t$ determine a potential
\begin{equation}\label{Ut}
\Ui_t(x) := \sum_{n\in\N} U_t(x-z^n_t), \qquad x\in\R^d.
\end{equation}
A point particle following the trajectory $f$ is then killed with rate $\Ui_t(f(t))$ at time $t$, and the probability that the particle has survived the traps by time $t$
is given by
$$
\exp\Big\{-\int_0^t \Ui_s(f(s))\, {\rm d}s\Big\}.
$$
Note that by replacing $U_t(\cdot)$ with $\widetilde U_t(\cdot):=U_t(\cdot+f(t))$, the problem is reduced to the case where the particle follows a constant trajectory. Therefore we may assume without loss of generality that $f\equiv 0$.

We are interested in upper boun{\rm d}s on the averaged survival probability
\be\label{St0}
S_t:= \E\Big[\exp\Big\{-\int_0^t \Ui_s(0) \, {\rm d}s\Big\}\Big]= \E\Big[\exp\Big\{-\sum_{n\in\N}\int_0^t U_s(-z^n_0-X^n_s)\, {\rm d}s\Big\}\Big],
\ee
where $\E$ denotes expectation w.r.t.\ $\Xi_0$ and $(X^n)_{n\in\N}$. Using Campbell's formula (see for instance~\cite[Section~3.2]{Kingman}) we obtain
\be\label{St}
S_t= \exp\Big\{ - \int_{\R^d} w_t(x) \phi(x) \, {\rm d}x \Big\},
\ee
where
\be\label{wt}
1-w_t(x):= \E_0\Big[\exp\Big\{ - \int_0^t U_s(-x-X_s)\, {\rm d}s\Big\} \Big] = \E_x\Big[\exp\Big\{ - \int_0^t U_s(-X_s)\, {\rm d}s\Big\} \Big],
\ee
where $\E_x$ denotes expectation w.r.t.\ the L\'evy process $X$ with $X_0=x$. We can interpret $w_t(x)$ as the probability that the L\'evy process $-X$, with $X_0=x$,
is killed before time $t$ by the trap $(U_t)_{t\geq 0}$. We will follow the convention that
\be\label{hardtrap}
\int_0^t U_s(-X_s) \, {\rm d}s:=\infty \qquad \mbox{if } U_s(-X_s)=\infty \mbox{ for some } s\in  [0,t),
\ee
so that the L\'evy process $-X_s$ is killed when it hits the hard trap $D_s:=\{x: U_s(x)=\infty\}$ for some $s<t$.

Analysis of the averaged survival probability $S_t$ then becomes equivalent to the analysis of
\be\label{Wt}
W^X_t(\phi, U_\cdot) := \int_{\R^d} w_t(x) \phi(x) \, {\rm d}x,
\ee
which can be interpreted as the total measure of $-X$ killed by the trap $U_\cdot:=(U_s)_{s\geq 0}$ up to time $t$, if $X$ starts with initial measure $\phi(x){\rm d}x$ on $\R^d$. Note that $W^X_t$ is exactly the continuous time analogue of $W_n$ in (\ref{ri}). In light of our discussion above, $e^{-W_n}$ can also be interpreted as
the averaged survival probability of a point particle in a Poisson field of moving traps in discrete time. As a corollary of Theorem~\ref{T:ri}, we will show that
$$
W^X_t(\phi, U_\cdot) \geq W^{X^*}_t(\phi_*, U^*_\cdot),
$$
where $U^*_\cdot:=(U^*_s)_{s\geq 0}$, and $X^*$ denotes the symmetric decreasing rearrangement of the L\'evy process $X$, which we now define.
\medskip

Recall that each L\'evy process $X$ with $X_0=0$ is uniquely characterized by a triple $(b, {\mathbb A}, \nu)$, called the characteristic of the L\'evy process (see e.g.~\cite{B96, S99}), such that the characteristic function of $X_t$ for any $t\geq 0$ is given by
$$
\E_0[e^{i \langle \xi, X_t\rangle}] = e^{-t\Psi(\xi)},
$$
where
\be\label{Psi}
\Psi(\xi) = -i \langle b, \xi\rangle + \frac{1}{2} \langle {\mathbb A}\xi, \xi\rangle + \int_{\R^d} (1-e^{i \langle \xi, x\rangle}+i \langle \xi, x\rangle 1_{\{|x|<1\}}) \nu({\rm d}x).
\ee
Here $b\in\R^d$ is a deterministic drift, $\mathbb A$ is the $d\times d$ covariance matrix of the Brownian component of $X$, and $\nu$ is a measure on $\R^d$ with
$$
\int_{\R^d} \frac{|x|^2}{1+|x|^2} \nu({\rm d}x) <\infty \qquad \mbox{and} \qquad \nu(\{0\})=0.
$$
The measure $\nu$ is called the L\'evy measure of $X$ and determines the jumps of $X$. When $b=0$, $\mathbb A=0$ and $\nu(\R^d)<\infty$, $X$ is simply
a compound Poisson process. Each L\'evy process admits a version with c\`adl\`ag sample paths, i.e., paths that are right continuous with left hand limits,
which we shall assume for $X$. If we denote by $\rho(x){\rm d}x$ the absolutely continuous part of $\nu$ with respect to the Lebesgue measure, then
the symmetric decreasing rearrangement of $X$ is defined to be the L\'evy process $X^*$ with characteristic $(0,{\mathbb A}^*, \nu^*)$, where ${\mathbb A}^*:={\rm Det}({\mathbb A})^{\frac{1}{d}}{\mathbb I}_d$ with ${\mathbb I}_d$ being the $d\times d$ identity matrix, and $\nu^*({\rm d}x)=\rho^*(x){\rm d}x$. This is the definition
of $X^*$ given in~\cite[Section 2]{W83}.

\brm\label{R:singular}{\rm
We note that the singular part of the L\'evy measure $\nu$ has been discarded in the definition of $X^*$. The reason is that we can rewrite the L\'evy process $X$ as the sum of two independent L\'evy processes $Y$ and $Z$, with $Z$ having the singular part of the L\'evy measure $\nu$. We can condition on $Z$ and treat it as a deterministic drift added to $Y$. Then the symmetric rearrangement of $Y+Z$ will lead to the removal of the drift $Z$, which is why we discard the singular part of $\nu$ (see Claim~\ref{clm:proofReduction} for more details). More generally, we can remove more than the singular part of the L\'evy measure before performing rearrangements. This will lead to a notion of domination for L\'evy processes.
}
\erm

We are now rea{\rm d}y to formulate our comparison result for the survival probability $S_t=e^{-W^X_t}$.
\bt\label{T:trap} Let $\phi: \R^d \to [0,\infty)$, and let $\phi_*$ be its symmetric increasing rearrangement defined in Theorem~\ref{T:ri}.
Let $U_\cdot(\cdot):[0,\infty) \times \R^d \to [0,\infty]$ be measurable, and for each $s\geq 0$, $|\{x: U_s(x)>l\}|<\infty$ for some $l<\infty$.
Assume that $D_s:=\{x: U_s(x)=\infty\}$ are open sets satisfying the regularity condition
$$
\mbox{\rm (R)} \quad \quad \forall\, s\geq 0 \mbox{ and } x\in D_s, \ \ \exists\, \delta>0, \ \ \mbox{s.t.} \ y \in D_{s'}\ \ \forall\, |y-x|<\delta
 \mbox{ and } s'\in [s, s+\delta). \qquad \qquad \quad
$$
Let $X$ be a L\'evy process with characteristic $(b, {\mathbb A}, \nu)$, and let $X^*$ be its symmetric decreasing rearrangement.
Let $W^X_t(\phi, U_\cdot)$ be defined from $X$, $\phi$, and $(U_s)_{s\geq 0}$ as in (\ref{Wt}), and let $W^{X^*}_t(\phi_*, U^*_\cdot)$ be defined analogously. Then for all $t\geq 0$,
\be\label{trap}
W^X_t(\phi, U_\cdot) \geq W^{X^*}_t(\phi_*, U^*_\cdot).
\ee
\et
\brm\label{R:trap}{\rm Condition (R) is equivalent to the lower semi-continuity of $f(s,x):=1_{\{x\in D_s\}}$ as $(t,y)\to (s^+, x)$. It guarantees that if the L\'evy process $-X_s \in D_s$ for some $s\in [0,t)$, then $-X_{s'}\in D_{s'}$ for all $s'\in [s, s+\delta)$
for some $\delta>0$, because $-X_s$ is almost surely right continuous in $s$. This ensures that our convention in (\ref{hardtrap}) is a.s.\ consistent with the usual definition of integral. The assumption on the level sets of $U_s$ will ensure that $U^*_\cdot$ also satisfies condition (R) (see the proof of (iii) of Claim~\ref{clm:proofReduction}).
Some natural sufficient conditions for (R) include: $D_s=D$ is an open set independent of time; $D_s=D+g(s)$ for an open set $D$ and a c\`adl\`ag path $g: [0,\infty)\to\R^d$; $\{(x,s) : s\geq 0, x\in D_s\}$ is an open set in $[0,\infty)\times\R^d$; $D_s^c$ is right continuous in $s$ with respect to the Hausdorff distance on the space of subsets of $\R^d$.
}\erm

The trapping problem defined above and its lattice version have been studied extensively in the physics literature, where the motion of the point particle can also be random (see e.g.~ \cite{BB02, MOBC04} and the references therein). It has also been studied as a detection problem in a mobile communication network (see e.g.~\cite{PSSS11} and the references therein); see also [CX11] for a recent study of the trapping problem with a renormalized Newtonian-type trap potential. A precursor to (\ref{trap}) in the literature is the special case when $X$ is a Brownian motion, $\phi\equiv 1$, and $U_s(x)=\infty \cdot 1_{\{|x+f(s)|< 1\}}$, where we recall that $f: [0,\infty)\to\R^d$ is the path of the point particle that was absorbed into the trap potential $(U_s)_{s\geq 0}$. In this case, inequality (\ref{trap}) only rearranges the function $f$. More precisely, it asserts that the survival probability $e^{-W^X_t}$ is maximized if the point particle follows the constant function $f\equiv 0$. This type of result, where the optimal trajectory is the constant trajectory, has been called the {\em Pascal principle} in the physics literature. For the lattice version of the trapping problem, the Pascal principle was established in~\cite{MOBC04}, see also~\cite[Corollary 2.1]{DGRS10}. In the continuum setting above where the spherical hard traps follow independent Brownian motions, it was first established in dimension 1 in~\cite{PSSS11}, assuming that $f$ is continuous. Subsequently, Peres and Sousi~\cite{PS11} generalized it to higher dimensions and proved (\ref{trap}) for the case where $X$ is a Brownian motion
and $U_s=\infty\cdot 1_{D_s}$ for arbitrary open sets $(D_s)_{s\geq 0}$. Their work and the work of Ba\~nuelos and M\'endez-Hern\'andez~\cite{BM-H10} motivated us to prove (\ref{trap}) in its current general form.
\medskip

Since the result of Peres and Sousi in~\cite{PS11} was formulated as an isoperimetric inequality for the expected volume of a Wiener sausage, which does not resemble (\ref{trap}) in appearance, we recall here the connection. In (\ref{trap}), let $\phi\equiv 1$ and let $U_s(\cdot):= \infty\cdot 1_{D_s}(\cdot)$, where $D_s:=D+g(s)$
for an open set $D\subset \R^d$ with finite volume and a c\`adl\`ag $g: [0,\infty)\to\R^d$. Note that $U_\cdot(\cdot)$ satisfies the assumptions in Theorem~\ref{T:trap}. From (\ref{wt}), we obtain
$$
w_t(x) = \P_x(-X_s\in D+g(s) \mbox{ for some } s\in [0,t)) = \P_0\Big(-x\in \bigcup_{s\in [0,t)} (D+X_s+g(s))\Big),
$$
where $\P_x$ denotes expectation for the L\'evy process $X$ with $X_0=x$. Then
$$
W^X_t(1, U_\cdot) = \!\int_{\R^d}\!\!\!\! w_t(x) \, {\rm d}x = \!\int_{\R^d} \!\!\! \E_0\big[1_{\{-x\in \bigcup_{s\in [0,t)} (D +X_s+g(s))\}}\big] \, {\rm d}x = \E_0\Big[{\rm Vol}\Big(\bigcup_{s\in [0,t)} \!\! (D +X_s+g(s)) \Big)\Big],
$$
where $\bigcup_{s\in [0,t)} (D+X_s+g(s))$ is the sausage generated by the L\'evy process $X$ with added drift $g$. Therefore in this case, (\ref{trap})
is equivalent to a comparison inequality for the expected volume of a L\'evy sausage. We formulate this as a Corollary.
\bcor\label{C:sausage}
Let $X$ be a L\'evy process and let $X^*$ be its symmetric decreasing rearrangement. Let $D\subset \R^d$ be an open set with finite volume, and let $g:[0,\infty)\to\R^d$
be c\`adl\`ag. Then for all $t>0$,
\be\label{sausage}
\E_0\Big[{\rm Vol}\Big(\bigcup_{s\in [0,t)} \!\! (D +X_s+g(s)) \Big)\Big] \geq \E_0\Big[{\rm Vol}\Big(\bigcup_{s\in [0,t)} \!\! (D^*+X^*_s) \Big)\Big].
\ee
\ecor
When $X$ is a Brownian motion, (\ref{sausage}) was proved in~\cite{PS11} with $D+g(s)$ replaced by any open set $D_s$, without even assuming the measurability of $D_s$ in $s$. We will not attempt such generality here due to the additional measure-theoretic complications it incurs.

\subsection{Comparison of Capacities}
As a corollary of Theorem~\ref{T:trap}, we establish a comparison inequality for the capacities of Borel sets for L\'evy processes. First we recall the definition of $q$-capacities for a L\'evy process $X:=(X_t)_{t\geq 0}$, with $q>0$ when $X$ is recurrent and $q\geq 0$ when $X$ is transient. Probabilistically, $q>0$ is the exponential rate of killing of the L\'evy process $X$, which ensures transience.

For any $A\in {\cal B}(\R^d)$, the Borel $\sigma$-algebra on $\R^d$, let $T_A(X):= \inf \{t \geq 0 : X_t \in A\}$ denote the first hitting time of $A$ by $X$. We will omit
$X$ from $T_A(X)$ when it is clear from the context with respect to which process $T_A$ is being evaluated. Let $\P_x(\cdot)$ and $\E_x[\cdot]$ denote probability and expectation for $X$ with $X_0=x$, and let $\widehat \P_x(\cdot)$ and $\widehat \E_x[\cdot]$ denote the analogues for $\widehat X:=-X$. We recall the following definition from~\cite[p.49]{B96} for $A$ either open or closed, and from \cite[Def.~6.1]{PS71} for general $A\in \Bi(\R^d)$.
\bdf\label{def:qCapa} {\bf ($q$-capacitary measure and $q$-capacities)}
Let $q>0$. For any $A\in \Bi(\R^d)$, the $q$-capacitary measure of $A$ for the L\'evy process $X$ is defined to be
\be\label{muqA}
\mu^q_A(B) := q \int_{\R^d} \E_x\big[e^{-qT_A} 1_{\{X_{T_A}\in B\}}\big] \, {\rm d}x \qquad \mbox{for all } B\in \Bi(\R^d).
\ee
Its total mass $C^q_X(A):=\mu^q_A(\R^d) = q\int_{\R^d} \E_x[e^{-qT_A}] \, {\rm d}x$ is called the $q$-capacity of $A$.
\edf
Some basic properties of $\mu^q_A$ and $C^q_X$ include:
\begin{itemize}
\item \cite[Thm.~6.2]{PS71} $\mu^q_A$ is the unique Radon measure supported on $\overline{A}$, the closure of $A$, with
\be\label{qpotential}
(\mu^q_A G^q)({\rm d}x) := \int_{\R^d} \mu^q_A({\rm d}y) G^q(y, {\rm d}x) = \widehat p^q_A(x) \, {\rm d}x,
\ee
where
\be\label{Gq}
G^q(y, B) := \int_0^\infty e^{-qt} \P_y(X_t\in B) \, {\rm d}t \qquad \mbox{for all } B\in \Bi(\R^d),
\ee
\be\label{pqA}
\widehat p^q_A(x) := \widehat \E_x[e^{-qT_A}].
\ee
Note that $G^q$ is the Green's function for the L\'evy process $X$, killed at exponential rate $q$.
\item \cite[Prop.~6.4]{PS71} $C^q_X(\cdot)$ is a Choquet capacity on $\Bi(\R^d)$. In particular, for all $A, B\in \Bi(\R^d)$,
\be\label{capincl}
C^q_X(A) \leq C^q_X(B) \qquad \mbox{if } A\subset B,
\ee
\be\label{qChoquet}
C^q_X(A) = \inf\{C^q_X(O): A\subset O, \, O\, \mbox{open}\} = \sup\{C^q_X(K) : K\subset A, \, K\, \mbox{compact}\}.
\ee
\end{itemize}
We note that~\eqref{qChoquet} is equivalent to $C^q_X$ being regular.

If $X$ is transient, i.e., $\lim_{t \to \infty} \vert X _t \vert = \infty$ a.s., then one can also define its $0$-capacity. We recall
the following definition from~\cite[Cor.~8, p.52]{B96} and \cite[Prop.~8.3]{PS71}.
\bdf\label{def:0Capa} {\bf ($0$-capacitary measure and $0$-capacities)} Suppose that $X$ is transient. Let $A\in \Bi(\R^d)$ be relatively
compact. Then $\mu^q_A$ converges weakly to a measure $\mu^0_A$, which is called the $0$-capacitary measure of $A$ for $X$. Its total mass
$C^0_X(A):=\mu^0_A(\R^d)$ is called the $0$-capacity (or just capacity) of $A$. For general $A\in \Bi(\R^d)$, we define
$C^0_X(A) := \sup\{ C^0_X(K) : K\subset A, K\, {\rm relatively\ compact}\}$.
\edf
For relatively compact $A\in\Bi(\R^d)$, the analogues of (\ref{qpotential}) and (\ref{capincl})--(\ref{qChoquet}) also hold~\cite[Prop.~8.2 \& 8.4]{PS71}, provided we replace $G^q$ by
\be\label{G0}
G^0(y, B) := \int_0^\infty \P_y(X_t\in B) \, {\rm d}t \qquad \mbox{for all } B\in \Bi(\R^d),
\ee
and replace $\widehat p^q_A(x)$ by
\be\label{p0A}
\widehat p^0_A(x) := \widehat \P_x(T_A<\infty).
\ee

We can now state our comparison inequality for capacities.
\bt\label{T:cap} Let $X$ be a L\'evy process with characteristic $(b, {\mathbb A}, \nu)$,
and let $X^*$ be its symmetric decreasing rearrangement.
Then for any $q>0$ ($q\geq 0$ if $X$ is transient), and for any $A \in \Bi(\R^d)$, we have
\be\label{cap}
C^q_X(A) \geq C^q_{X^*}(A^*).
\ee
\et
\brm{\rm
Theorem~\ref{T:cap} was conjectured in \cite[p.4050]{BM-H10}. It exten{\rm d}s a result of Watanabe~\cite[Theorem~1]{W83}, where (\ref{cap})
was proved for symmetric L\'evy processes, i.e., $X$ is equally distributed with $-X$ if $X_0=0$. Classic Newtonian capacities correspond to $X$ being a Brownian motion.
For Riesz capacities, which correspond to radially symmetric $\alpha$-stable processes, Watanabe's result has been rediscovered by Betsakos in~\cite{B04}. For isotropic unimodal L\'evy processes, Watanabe's result has been rediscovered by M\'endez-Hern\'andez in~\cite{M-H06}, which uses the Friedberg-Luttinger inequality
discussed in Remark~\ref{R:FL}.}
\erm

In~\cite{W83}, Watanabe used the definition of $q$-capacities from the theory of Dirichlet forms for symmetric Markov processes. It is known that such a definition is equivalent to the probabilistic definition given here if $X$ is a symmetric L\'evy process. However, a precise reference seems hard to locate. Therefore we will sketch briefly why the two definitions are equivalent.

For a symmetric L\'evy process $X$ with characteristic $(0, {\mathbb A}, \nu)$, one can define a family of
Dirichlet forms $\Ei_q(\cdot, \cdot)$ ($q>0$ if $X$ is recurrent and $q\geq 0$ if $X$ is transient). If $L$ denotes the generator
of $X$, and $L_q u:= L u -qu$ for $u\in \Di(L)\subset L^2(\R^d)$, then the domain of $\Ei_q$ equals $\Di(\Ei_q)=\Di(\sqrt{-L_q})\subset L^2(\R^d)$, and
$$
\begin{aligned}
\Ei_q(u, v) & = \int_{\R^d} (\sqrt{-L_q} u)(x) (\sqrt{-L_q} v)(x) \, {\rm d}x  \qquad && \mbox{for } u,v\in \Di(\Ei_q), \\
\Ei_q(u, v) & = - \int_{\R^d} u(x) (L_qv)(x) \, {\rm d}x  \qquad && \mbox{for } u\in \Di(\Ei_q), v\in \Di(L_q);
\end{aligned}
$$
see Theorem~1.3.1 and Corollary~1.3.1 in \cite{FOT11}. For any open set $O$, its $q$-capacity is defined by (see
e.g.~\cite{W83} or~\cite[Chap.~2]{FOT11})
\be\label{Ecap}
C^q_X(O):= \inf\{\Ei_q(u, u) : u\in \Di(\Ei_q), u\geq 1 \mbox{ a.e.~on } O\},
\ee
with $C^q_X(O):=\infty$ if the infimum is taken over an empty set. The $q$-capacity of a general Borel set $A\subset\R^d$ is defined to be
\be\label{Ecap2}
C^q_X(A) := \inf\{ C^q_X(O) : A\subset O, O \mbox{ open}\}.
\ee
For any bounded open set $O$, Lemma~2.1.1 (and the remark before Lemma~2.1.8) in \cite{FOT11} show that the infimum in (\ref{Ecap}) is achieved at a function $e^q_O\in \Di(\Ei_q)$, called the $q$-equilibrium potential of $O$. Furthermore,
\be\label{Ecap3}
C^q_X(O) = \Ei_q(e^q_O, v) = \int_{\R^d} e^q_O(x) (-L_qv)(x) \, {\rm d}x \qquad \forall\, v\in \Di(L_q) \mbox{ with } v= 1 \mbox{ a.e. on } O.
\ee
Lemma 4.2.1 and Theorem 4.3.3 in~\cite{FOT11} identify $e^q_O$ with $\widehat p^q_O$ defined in (\ref{pqA}) and (\ref{p0A}). Therefore, choosing any $v\in C^\infty_c(\R^d)\subset \Di(L_q)$ with $v=1$ on $\overline{O}$, we obtain
\be\label{Ecap4}
C^q_X(O) = \int \widehat p^q_O(x) (-L_qv)(x) \, {\rm d}x = \int \mu^q_O({\rm d}y) \int G^q(y, {\rm d}x) (-L_q v)(x) = \int \mu^q_O({\rm d}y) v(y),
\ee
where we used (\ref{qpotential}), and the fact that $G^q(x, {\rm d}y)$ is the Green's kernel for the L\'evy process $X$ killed with rate $q$, which is a transient process with generator $L_q$, and hence $G^q(x, {\rm d}y)$ defines an integral operator which is the inverse of $-L_q$. Since $v=1$ on $\overline O$ and $\mu^q_O$ is supported
on $\overline{O}$, $C^q_X(O)$ in (\ref{Ecap4}) coincides with our definition of $C_X^q(O)$ in Definitions~\ref{def:qCapa}--\ref{def:0Capa}.
Since $C^q_X(\cdot)$ defined via (\ref{Ecap})--(\ref{Ecap2}) is also a Choquet capacity by~\cite[Theorem 2.1.1]{FOT11} and hence satisfies (\ref{qChoquet}), the coincidence of the two definitions of capacities extends from bounded open sets to all Borel-measurable sets.

\subsection{Outline}

The rest of the paper is organized as follows. Theorems~\ref{T:ri}, \ref{T:trap} and \ref{T:cap} will be proved respectively in Sections~\ref{S:ri}, \ref{S:trap} and \ref{S:cap}. Section~\ref{sec:conclusion} discusses some open questions. Lastly, in Appendix~\ref{S:conv}, we collect some basic properties of symmetric rearrangements that we use in the proofs.

\section{Proof of Theorem~\ref{T:ri}}\label{S:ri}

The proof of Theorem~\ref{T:ri} is surprisingly simple. With the introduction of a suitable notion of symmetric domination (see (\ref{dom})) on the density profile of the surviving random walk, Theorem~\ref{T:ri} is then deduced by induction over the number of time steps. Lemmas~\ref{L:domV} and \ref{L:domp} are the key lemmas, and the only tool we need to use for their proofs is the Riesz rearrangement inequality.
\bigskip

\noindent
{\bf Proof of Theorem~\ref{T:ri}.} By scaling and replacing $\phi$ with $\phi\wedge\sigma$, we may assume without loss of generality that $\sigma=1$, $\phi \in [0,1]$,
and $|\{x: \phi(x)<t\}|<\infty$ for all $t\in [0,1)$. By truncating $(V_i)_{i\geq 0}$ and $(p_i)_{i\geq 1}$ and then applying Lemma~\ref{L:approx} and the Monotone Convergence Theorem, we may first assume without loss of generality that $(V_i)_{i\geq 0}$ and $(p_i)_{i\geq 1}$ are integrable. Furthermore, we may assume that $(p_i)_{i\geq 1}$ are probability densities. For such $(V_i)_{i\geq 0}$ and $(p_i)_{i\geq 1}$, we can then apply Lemma~\ref{L:approx} and the Dominated Convergence Theorem to reduce to the case where $1-\phi$ is integrable, which we assume from now on.

Let us denote $\psi:=1-\phi$. Since $(p_i)_{i\geq 1}$ are assumed to be probability densities, we can rewrite $W_n(\phi, V_\cdot, p_\cdot)$ in (\ref{ri}) as
\begin{eqnarray}
W_n(\phi, V_\cdot, p_\cdot)\!\!\!\! & = &\!\!\!\!\! \int\!\!\cdots\!\!\int  (1-\psi(x_0)) \Big\{1-\prod_{i=0}^n (1-V_i(x_i))\Big\}\prod_{i=1}^n p_i(x_i-x_{i-1})   \prod_{i=0}^n {\rm d}x_i \label{psisplit}\\
\nn \!\!\!&=&\!\!\!\! - \!\! \int\!\!\psi(x_0)\, {\rm d}x_0 + \!\int\!\!\cdots\!\!\int\!\! \Big\{1-(1-\psi(x_0))\prod_{i=0}^n (1-V_i(x_i))\Big\}\!\prod_{i=1}^n p_i(x_i-x_{i-1})   \prod_{i=0}^n {\rm d}x_i.
\end{eqnarray}
A similar identity hol{\rm d}s for $W_n(\phi_*, V^*_\cdot, p^*_\cdot)$, since $(p_i^*)_{i\geq 1}$ are also probability densities.
We will be guided by the probabilistic interpretation that $(p_i)_{i\geq 1}$ are the transition probability densities of a random walk $X$, which
is killed at each time $i\geq 0$ with probability $V_i(X_i)$. We can also interpret $\psi$ as a trap function at time $0$, so that $X$ is killed at time $0$ first with probability $\psi(X_0)$, and in case it survives, it is then killed with probability $V_0(X_0)$. If we start $X$ at time $0$ with Lebesgue measure, then
$$
\phi_0(x_0):= (1-V_0(x_0))(1-\psi(x_0))
$$
is the density of $X_0$ on $\R^d$ upon surviving the traps $\psi$ and $V_0$. Similarly, for $n\in\N$,
\be\label{phin}
\phi_n(x_n) := (1-V_n(x_n))\idotsint (1-\psi(x_0))\prod_{i=0}^{n-1} (1-V_i(x_i))\prod_{i=1}^n p_i(x_i-x_{i-1}) \prod_{i=0}^{n-1} {\rm d}x_i
\ee
is the density of $X_n$ on $\R^d$ upon survival up to time $n$. We can then rewrite (\ref{psisplit}) as
$$
W_n(\phi, V_\cdot, p_\cdot) + \int\psi(x_0)\, {\rm d}x_0 = \int (1-\phi_n(x_n)) \,{\rm d}x_n,
$$
which is the total measure of $X$ killed up to time $n$. Similarly,
$$
W_n(\phi_*, V^*_\cdot, p^*_\cdot) + \int\psi^*(x_0) \, {\rm d}x_0 = \int (1-\vartheta_n(x_n)) \, {\rm d}x_n,
$$
where
\be\label{varphin}
\vartheta_n(x_n) := (1-V^*_n(x_n))\idotsint (1-\psi^*(x_0))\prod_{i=0}^{n-1} (1-V^*_i(x_i))\prod_{i=1}^n p^*_i(x_i-x_{i-1}) \prod_{i=0}^{n-1} \, {\rm d}x_i.
\ee
Since $\int\psi = \int\psi^*$, to prove (\ref{ri}), it then suffices to show that
\be\label{riphi}
\int (1-\vartheta_n(x_n)) \, {\rm d}x_n \leq \int (1-\phi_n(x_n)) \, {\rm d}x_n.
\ee
The key insight in our proof (\ref{riphi}) is the introduction of a notion of {\em symmetric domination}. More precisely, we show that $\vartheta_n$ symmetrically dominates $\phi_n$, denoted by $\vartheta_n\succ\phi_n$, in the sense that
\be\label{dom}
\int_{(A^*)^c} (1-\vartheta_n(x_n)) \, {\rm d}x_n \leq \int_{A^c} (1-\phi_n(x_n)) \, {\rm d}x_n \qquad \mbox{for all measurable } A \mbox{ with } |A|<\infty.
\ee
Heuristically, this means that $\vartheta_n(x){\rm d}x$ contains more mass and with mass closer to $\infty$ than the symmetric decreasing rearrangement of $\phi_n(x){\rm d}x$ around $\infty$. Note that $\succ$ is a partial order on the class of functions $f: \R^d\to [0,1]$ with $\int (1-f)<\infty$. By setting $A=\{0\}$ in (\ref{dom}), we obtain (\ref{riphi}). We remark that closely related notions of symmetric domination have appeared in various contexts, such as in~\cite{ATL89} for the study of rearrangement inequalities around $0$ instead of $\infty$, or in~\cite{K77} for the study of symmetric unimodal distributions and Anderson's inequality, which involves comparison of measures on the complements of symmetric convex sets, but without rearrangements.

Note that
for $n\in\N$,
\be\label{phirec}
\begin{aligned}
\phi_n(x_n) & = (1-V_n(x_n)) (p_n*\phi_{n-1})(x_n), \\
\vartheta_n(x_n) & = (1-V^*_n(x_n)) (p^*_n*\vartheta_{n-1})(x_n).
\end{aligned}
\ee
Therefore by induction, to prove $\vartheta_n\succ\phi_n$, it suffices to show that: $(1-\psi^*)\succ (1-\psi)$ (recall that $\psi:=1-\phi$); and if $\vartheta\succ\phi$, then
$(1-V^*)\vartheta \succ (1-V)\phi$ and $p^**\vartheta \succ p*\phi$ for any integrable $V:\R^d\to [0,1]$ and any probability density $p:\R^d\to [0,\infty)$.
The first fact holds because for any measurable $A$ with $|A|<\infty$,
$$
\int_{A^c} \psi(x) \, {\rm d}x = \int \psi(x) \, {\rm d}x - \int 1_{A}(x) \psi(x) \, {\rm d}x \geq \int \psi^*(x) \, {\rm d}x - \int 1_{A^*}(x) \psi^*(x) \, {\rm d}x = \int_{(A^*)^c} \psi^*(x)\, {\rm d}x
$$
by a classic rearrangement inequality (see e.g.~\cite[Theorem 3.4]{LL01}). The other claims on the preservation of $\succ$ hold by Lemmas~\ref{L:domV} and \ref{L:domp} below, where the integrability conditions therein are guaranteed by our integrability assumptions on $\psi$ and $(V_i)_{i\geq 0}$.
\qed

We now state and prove the two key lemmas used in the proof of Theorem~\ref{T:ri}.

\blem\label{L:domV} Suppose that $\phi, \vartheta, V: \R^d\to [0,1]$ are such that $(1-\phi)$, $(1-\vartheta)$ and $V$ are all integrable. If $\vartheta\succ \phi$ in the sense defined in (\ref{dom}), then we also have $(1-V^*)\vartheta \succ (1-V)\phi$.
\elem
{\bf Proof.} We need to show that for all measurable $A$ with $|A|<\infty$,
\be\label{domv1}
\int_{(A^*)^c} (1-(1-V^*)\vartheta) \leq \int_{A^c} (1-(1-V)\phi).
\ee
By writing $V(x)=\int_0^1 1_{F_t}(x) \, {\rm d}t$ with $F_t:=\{x: V(x)>t\}$, and $V^*(x) = \int_0^1 1_{F^*_t}(x) \, {\rm d}t$, where we note that $F^*_t:=\{x: V^*(x)>t\}$
is also the symmetric decreasing rearrangement of $F_t$, it suffices to verify (\ref{domv1}) for the case $V=1_F$ for some measurable set $F$ with $|F|<\infty$.
For $V=1_F$, the LHS of (\ref{domv1}) equals
\be\label{domv2}
\int_{(A^*)^c} (1-1_{(F^*)^c}\vartheta) = \int_{(A^*)^c} (1_{F^*}+1_{(F^*)^c}(1-\vartheta)) = |F^*|-|F^*\cap A^*| +\int_{(F^*\cup A^*)^c} (1-\vartheta).
\ee
Similarly, the RHS of (\ref{domv1}) equals
\be\label{domv3}
\int_{A^c} (1-1_{F^c} \phi) = |F|-|F\cap A| +\int_{(F\cup A)^c} (1-\phi).
\ee
Since $|F^*|=|F|$, and the remaining terms in (\ref{domv2}) and (\ref{domv3}) are symmetric in $F$ and $A$, we may assume without loss of generality that
$|F|\geq |A|$, which implies $A^*\subset F^*$. Subtracting (\ref{domv2}) from (\ref{domv3}) then gives
\be
\begin{aligned}
& -|F\cap A| + |A^*| + \int_{(F\cup A)^c} (1-\phi) - \int_{(F^*)^c} (1-\vartheta) \\
=\ & |F^c\cap A| - \int_{F^c\cap A} (1-\phi) + \int_{F^c}(1-\phi) - \int_{(F^*)^c} (1-\vartheta) \geq\ \int_{F^c\cap A} \phi \geq 0,
\end{aligned}
\ee
where we used the fact that $|A^*|=|A|$, and the assumption $\vartheta\succ \phi$. This proves (\ref{domv1}).
\qed

\blem\label{L:domp} Suppose that $\phi, \vartheta: \R^d\to [0,1]$ are such that $(1-\phi)$ and $(1-\vartheta)$ are integrable. If $\vartheta\succ \phi$, then for any probability
density $p: \R^d \to [0,\infty)$, we have $p^* * \vartheta \succ p*\phi$.
\elem
{\bf Proof.} First we note that $\vartheta \succ \phi_* \succ \phi$, where $\phi_*:=1-(1-\phi)^*$ is the symmetric increasing rearrangement of $\phi$. This follows
from the observation that $\vartheta \succ \phi$ implies
$$
\int_{(A^*)^c} (1-\vartheta) \leq \inf_{B: |B|=|A|} \int_{B^c} (1-\phi) = \int_{(A^*)^c} (1-\phi)^* = \int_{(A^*)^c} (1-\phi_*) \leq \int_{A^c} (1-\phi_*).
$$
For any measurable set $A$ with $|A|<\infty$, we have
\begin{eqnarray*}
\int_{A^c} (1-p*\phi) = \int_{A^c} p*(1-\phi) &=& \int(1-\phi) - \iint 1_A(x) p(x-y)(1-\phi)(y) \, {\rm d}y \, {\rm d}x \\
&\geq&  \int (1-\phi)^* - \iint 1_{A^*}(x) p^*(x-y)(1-\phi)^*(y) \, {\rm d}y \,{\rm d}x \\
&=& \int (1-\phi)^*(1- p^**1_{A^*}),
\end{eqnarray*}
where in the inequality we used Riesz's rearrangement inequality~\cite[Theorem 3.7]{LL01}. Note that
$$
1-p^* * 1_{A^*}(y) = 1-\int_0^1 1_{\{p^**1_{A^*}(y)>t\}} \, {\rm d}t = \int_0^1 1_{\{p^**1_{A^*}(y)\leq t\}} \, {\rm d}t,
$$
where $\{y: p^**1_{A^*}(y)>t\}$, $t\in (0,1)$, are centered open balls because $p^**1_{A^*} = (p^**1_{A^*})^*$ by Lemma~\ref{L:conv}. Since $\vartheta\succ \phi_*$,
we then have
\begin{eqnarray*}
\int_{A^c} (1-p* \phi) \geq \int (1-\phi)^*(1- p^**1_{A^*}) &=& \int_0^1 \int_{\{y: p^**1_{A^*}(y)\leq t\}} (1-\phi)^*(y) \, {\rm d}y\, {\rm d}t \\
&\geq& \int_0^1 \int_{\{y: p^**1_{A^*}(y)\leq t\}} (1-\vartheta)(y) \, {\rm d}y\, {\rm d}t \\
&=& \int (1-\vartheta) (1- p^**1_{A^*}) = \int (1-\vartheta) p^**1_{(A^*)^c} \\
&=& \int 1_{(A^*)^c} p^**(1-\vartheta) = \int_{(A^*)^c} (1-p^**\vartheta).
\end{eqnarray*}
Therefore $p^* * \vartheta \succ p * \phi$.
\qed

\brm\label{R:riext}{\rm Theorem~\ref{T:ri} admits two extensions which follow by the same proof as above. Firstly, (\ref{ri}) remains valid if for each $i\geq 0$, we replace $(1-V_i(x_i))$ by $\prod_{k=1}^{l_i} (1-V^{(k)}_i(x_i))$ for some $l_i\in\N$ and $V^{(k)}_i: \R^d\to [0,1]$ for $1\leq k\leq l_i$, and replace $(1-V^*_i(x_i))$ by $\prod_{k=1}^{l_i} (1-V^{(k)*}_i(x_i))$. Secondly, assuming $\sigma= 1$ in Theorem~\ref{T:ri} and $\int (1-\phi)<\infty$, then (\ref{ri}) also holds if we replace $\phi_*$ by any $\vartheta: \R^d\to [0,1]$ such that $\int (1-\vartheta)=\int(1-\phi)$, and $\vartheta$ symmetrically dominates $\phi$ in the sense defined in (\ref{dom}). This latter extension also applies to Theorem~\ref{T:trap}.
}
\erm

\section{Proof of Theorem~\ref{T:trap}}\label{S:trap}

We will first make the reductions stated below in Claim~\ref{clm:proofReduction}, and then perform discrete time approximation and apply Theorem~\ref{T:ri}.
\begin{clm} \label{clm:proofReduction}
It is sufficient to prove Theorem~\ref{T:trap} for
\begin{enumerate} [label=(\roman{*}), ref=(\roman{*})]
\item[\rm (i)]
$\phi \in [0,1]$,

\item[\rm (ii)]
L\'evy measures $\nu({\rm d}x)=\rho(x){\rm d}x$ for some $\rho: \R^d\to [0,\infty)$,

\item[\rm (iii)] potentials
$U: [0,\infty) \times \R^d \to [0,\infty)$ which are continuous with bounded support.
\end{enumerate}
\end{clm}
{\bf Proof of (i).} This follows by the same reasoning as in the proof of Theorem~\ref{T:ri}. \qed
\medskip

\noindent
{\bf Proof of (ii).} Assume that Theorem~\ref{T:trap} holds under assumption (ii). As noted in Remark~\ref{R:singular}, for a general L\'evy process $X$ with characteristic $(b, \mathbb A, \nu)$, with $\rho(x){\rm d}x$ being the absolutely continuous part of $\nu$, $X$ is equally distributed with $Y+Z$, where $Y$ is a L\'evy process with characteristic $(b, \mathbb A, \rho(x) {\rm d}x)$ and $Y_0=X_0$, and $Z$ is an independent L\'evy process with characteristic $(0,0, \nu - \rho(x) {\rm d}x)$ and $Z_0=0$. Let $\E^Y_y$ denote
expectation for $Y$ with $Y_0=y$, and let $\E^Z_0$ be defined similarly. By Tonelli's Theorem, we have
\begin{equation}\label{3.1}
\begin{aligned}
W_t^X (\phi, U_\cdot)
&= \int_{\R^d} \phi(x)  \Big(1- \E_x \Big[ e^{ - \int_0^t U_s(-X_s) \, {\rm d}s} \Big] \Big) \, {\rm d}x\\
&= \E_0^Z \Bigg[ \int_{\R^d} \phi(x)  \Big(1- \E^Y_x \Big[ e^{ - \int_0^t U_s(-Y_s-Z_s) \, {\rm d}s }  \Big] \Big) \, {\rm d}x \Bigg]\\
&\ge \int_{\R^d} \phi_*(x)  \Big (1- \E^{Y^*}_x \Big[ e^{ - \int_0^t U^*_s(-Y^*_s) \, {\rm d}s } \Big] \Big) \, {\rm d}x
= W_t^{X^*} (\phi_*,U_\cdot^*).
\end{aligned}
\end{equation}
In the inequality above, conditional on $Z$, we applied Theorem~\ref{T:trap} for the L\'evy process $Y$ which satisfies assumption (ii), and we applied
symmetric decreasing rearrangement to the potential $\widetilde U_s(x) := U_s(x-Z_s)$. Note that because $Z$ is a.s.\ c\`adl\`ag, $\widetilde U_\cdot$ also satisfies
the regularity condition (R) in Theorem~\ref{T:trap}; furthermore, we note that $\widetilde U^*_\cdot = U^*_\cdot$. In the last equality above, we used the fact
that $Y^*$ and $X^*$ are equal in law. This proves the reduction to L\'evy processes satisfying (ii).
\qed
\medskip

\noindent
{\bf Proof of (iii).} We first reduce to potentials $U:[0,\infty)\times\R^d\to [0,\infty)$ which are bounded with bounded support. Assume that Theorem~\ref{T:trap} holds
for such potentials. For a general potential $U$ satisfying the conditions in Theorem~\ref{T:trap}, and for each $n\in\N$, define $U_{n,s}(x) := 1_{\{s+|x|<n\}} U_s(x) \wedge n$. Then $U_{n,\cdot}$ is bounded with bounded support, and $U_{n,s}(x)\uparrow U_s(x)$ for all $s\geq 0$ and $x\in\R^d$ as $n\uparrow\infty$.

If Theorem~\ref{T:trap} holds for bounded potentials with bounded support, then
\be\label{WUncomp}
W^X_t(\phi, U_{n, \cdot}) \geq W^{X^*}_t (\phi_*, U^*_{n, \cdot}) \qquad \mbox{for all } n\in\N.
\ee
To prove $W^X_t(\phi, U_\cdot)\geq W^{X^*}_t(\phi_*, U^*_\cdot)$, and thus
complete the reduction to bounded potentials with bounded support, it suffices to show that
\be\label{WUnmono}
W^X_t(\phi, U_{n,\cdot}) \ \big\uparrow \ W^X_t(\phi, U_\cdot) \qquad \mbox{as } n\uparrow \infty
\ee
and
\be\label{WUn*mono}
W^{X^*}_t(\phi_*, U^*_{n,\cdot}) \ \big\uparrow \ W^{X^*}_t(\phi_*, U^*_\cdot) \qquad \mbox{as } n\uparrow \infty.
\ee
We first claim that for every $x\in\R^d$ and for almost every realization of $X$ with $X_0=x$,
\be\label{Unmono}
\int_0^t U_{n,s}(-X_s) \, {\rm d}s\ \Big\uparrow\ \int_0^t U_s(-X_s) \, {\rm d}s \qquad \mbox{as } n\uparrow \infty.
\ee
Indeed, if $U_s(-X_s)<\infty$ for all $s\in [0,t)$, then (\ref{Unmono}) follows by the Monotone Convergence Theorem; if $U_s(-X_s)=\infty$ for some $s\in [0,t)$,
so that $\int_0^t U_s(-X_s) \, {\rm d}s :=\infty$ by our convention in (\ref{hardtrap}), then the regularity assumption (R) in Theorem~\ref{T:trap} and Remark~\ref{R:trap} imply that (\ref{Unmono}) still holds. Applying the Monotone Convergence Theorem to the expression for $W^X_t$, c.f.~(\ref{3.1}), then gives (\ref{WUnmono}).

To verify (\ref{WUn*mono}), we note that by Lemma~\ref{L:approx}, we have $U^*_{n,s}(x)\uparrow U^*_s(x)$ for all $s\geq 0$ and $x\in\R^d$ as $n\uparrow \infty$. Furthermore, the potential $U^*_\cdot$ also satisfies condition (R) in Theorem~\ref{T:trap}. This is because $U_\cdot$ satisfies (R), which implies that its infinity level sets $(D_s)_{s\geq 0}$ satisfy
$$
1_{D_s}(x) \leq \liminf_{s'\downarrow s} 1_{D_{s'}}(x) \qquad \mbox{for all } x\in\R^d, s\geq 0.
$$
Therefore by Fatou's lemma, $|D_s| \leq \liminf_{s'\downarrow s} |D_{s'}|$ for all $s\geq 0$.
 The assumption in Theorem~\ref{T:trap} that $|\{x: U_s(x)>l\}|<\infty$
for some $l<\infty$ implies that $|D_s| < \infty$ and $\{x: U^*_s(x)=\infty\}=D_s^*$, so $|D^*_s| \leq \liminf_{s'\downarrow s} |D^*_{s'}|$. Since $(D_s^*)_{s\geq 0}$ are finite centered open balls, it is easily seen that $U^*_\cdot$ must also satisfy condition (R). The same arguments as those leading to (\ref{WUnmono}) then imply (\ref{WUn*mono}).

We now make the further reduction from bounded $U$ with bounded support to continuous $U$ with bounded support. For any bounded $U_\cdot(\cdot)$ with bounded support, we can find a sequence of continuous $U_{n,\cdot}(\cdot)$, uniformly bounded with uniformly bounded support, such that for all $(s,x)$ in a set $N\subset [0,\infty)\times\R^d$ with
full Lebesgue measure, we have $U_{n,s}(x)\to U_s(x)$ as $n\to\infty$. By Fubini's Theorem, for every realization of the L\'evy process $X$ with $X_0=0$, we have
$$
0=\int_0^\infty \int_{\R^d} 1_{N^c}(s,x) \, {\rm d}x \, {\rm d}s = \int_0^\infty \int_{\R^d} 1_{N^c}(s,-x-X_s) \, {\rm d}x\, {\rm d}s = \int_{\R^d} \int_0^\infty  1_{N^c}(s,-x-X_s) \, {\rm d}s\, {\rm d}x.
$$
Therefore for every $x$ in a set $\Lambda\subset\R^d$ with $|\Lambda^c|=0$, the set $\{s\geq 0: (s, -x-X_s)\in N\}$ has full Lebesgue measure on $[0,\infty)$. Writing a L\'evy process starting from $x$ as $x$ plus a L\'evy process starting from the origin, we can write
\be\label{domconv}
\begin{aligned}
W_t^X (\phi, U_{n,\cdot})
=&\ \int_{\Lambda} \phi(x) \Big( 1- \E_x \Big[ e^{ - \int_0^t U_{n,s} (-X_s) \, {\rm d}s } \Big] \Big) \,{\rm d}x \\
=&\ \E_0 \Big[ \int_{\Lambda} \phi(x) \Big( 1- e^{ - \int_0^t U_{n,s} (-x-X_s) 1_{N} (s,-x-X_s) \, {\rm d}s} \Big) \,{\rm d}x  \Big] \\
\asto{n} &\ \E_0 \Big[ \int_{\Lambda} \phi(x) \Big( 1- e^{ - \int_0^t U_s(-x-X_s) 1_{N} (s,-x-X_s)\, {\rm d}s } \Big) \,{\rm d}x  \Big] \\
=&\ \int_{\Lambda} \phi(x) \Big( 1- \E_x \Big[ e^{ - \int_0^t U_s (-X_s)  \, {\rm d}s}  \Big] \Big) \,{\rm d}x = W_t^X (\phi, U_{\cdot}).
\end{aligned}
\ee
The convergence above holds by the Dominated Convergence Theorem because the integrands under $\E_0\big[\int_{\Lambda}\cdot \big]$ can be dominated uniformly by
$1_{\{T_B(-x-X)< t\}}$, where $B$ is a finite open ball containing the support of $U_{n,s}$ and $U_s$ for all $s\geq 0$ and $n\in\N.$ Note that
$$
\E_0\Big[\int_{\Lambda}1_{\{T_B(-x-X)< t\}} \, {\rm d}x \Big] = \int_{\R^d} \P_x(T_B(-X)<t) \, {\rm d}x,
$$
which is finite by~\cite[Prop.~3.6]{PS71}, and hence the Dominated Convergence Theorem can be applied.

By Lemma~\ref{lem:conv}, for Lebesgue a.e.\ $s\geq 0$, we have $U^*_{n,s}(x)\to U^*_s(x)$ for Lebesgue a.e.\ $x\in\R^d$. Therefore we can apply the same argument as above to conclude that $W_t^{X^*} (\phi_*, U^*_{n,\cdot})\to W_t^{X^*} (\phi_*, U^*_{\cdot})$ as $n\to\infty$. If Theorem~\ref{T:trap} holds for continuous potentials with bounded support, then we have $W_t^X (\phi, U_{n,\cdot})\geq W_t^{X^*} (\phi_*, U^*_{n,\cdot})$, which as $n\to\infty$ implies the same comparison for $U$. Therefore Theorem~\ref{T:trap} also holds for bounded potentials with bounded support, which concludes the reduction to potentials satisfying (iii).
\qed
\medskip

To prove Theorem~\ref{T:trap} under the assumptions in Claim~\ref{clm:proofReduction}, we will follow the same steps as in~\cite{BM-H10}. We will discretize time, and
approximate the L\'evy process $X$ with characteristic $(b, {\mathbb A}, \rho(x){\rm d}x)$ in the standard way by truncating its L\'evy measure $\rho(x){\rm d}x$, so that we have the sum of a compound Poisson process and an independent Brownian motion. For this purpose we define $\rho_n(y) := \rho(y) 1_{ \{ |y|>1/n \} }$ and let $c_n:= \int_{\R^d}\rho_n(y) \, {\rm d}y$, so that $\bar\rho_n(y):=c_n^{-1}\rho(y)$ is a probability density on $\R^d$. Let $C_{n,t}$ be a compound Poisson process, starting at $0$, with characteristic function
\[
\E_0[e^{i\langle \xi, C_{n,t}\rangle}] = e^{-t \bar\Psi_n(\xi)},
\]
where
\[
\bar\Psi_{n}(\xi) = \int_{\R^d} (1-e^{i\langle \xi, y\rangle}) \rho_n(y) \, {\rm d}y = c_n \int_{\R^d} (1-e^{i\langle \xi, y\rangle}) \bar\rho_n(y) \, {\rm d}y.
\]
Choose $\epsilon_n$ to be a sequence of positive numbers converging to $0.$ Then with ${\mathbb I}_d$ denoting the $d\times d$ identity matrix, ${\mathbb A}_n := {\mathbb A}+ \epsilon_n {\mathbb I}_d$ is a positive definite matrix since $\mathbb A$ is positive semi-definite. Let $G_{n,t}$ be a
Brownian motion independent of $C_{n,t}$, starting at $x$, with covariance matrix ${\mathbb A}_n$ and drift $b_n=b-\int_{|y|<1} y \rho_n(y) \, {\rm d}y$. Now set $X_{n,t}:=C_{n,t} + G_{n,t}$.
Since $C_{n,t}$ and $G_{n,t}$ are independent, we get that
\[
\E_x[e^{i\langle \xi, X_{n,t}\rangle}] = e^{-t \Psi_n(\xi) + i \langle\xi,x\rangle},
\]
where
\[
\Psi_n(\xi) = -i \langle b, \xi \rangle + \frac{1}{2}\langle {\mathbb A}_n \xi, \xi \rangle
+ \int_{\R^d} \big( 1+i\langle \xi,y \rangle 1_{\{ |y|<1\}} - e^{i\langle \xi, y\rangle} \big) \rho_n(y) \, {\rm d}y.
\]
We first prove a discrete time analogue of Theorem~\ref{T:trap} for $X_{n,\cdot}:=(X_{n,t})_{t\geq 0}$, which approximates $X$.
\blem\label{lem:gaussian}
Let $X_{n,\cdot}$ be as above. Let $\phi: \R^d\to [0,\infty)$, and let $m\in\N$. For $1\leq i\leq m$, let $V_i: \R^d\to [0,1]$ be continuous with compact support. Let
$0 < t_1 < \ldots < t_m<\infty$. Then
\begin{align}\label{eq:approx}
\int_{\R^d} \phi(x) \Big(1- \E_x\Big[\prod_{i=1}^{m}(1-V_i(X_{n,t_i})) \Big] \Big) \, {\rm d}x
\geq \int_{\R^d} \phi_*(x) \Big(1- \E_x \Big[ \prod_{i=1}^{m}(1-V^*_i(X^*_{n,t_i})) \Big] \Big)\, {\rm d}x.
\end{align}
\elem
{\bf Proof.} Let $p_{n,t}(\cdot)$ denote the transition kernel of $G_{n,t}$. With the convention that $t_0 =0$, $k_0=0$, $z_0 = x,$ and denoting by $N_{n,t}$ the Poisson process
which counts the number of jumps of $C_{n,t}$, we can write
\begin{align} \label{eq:summandRep}
\begin{split}
1-\E_x \Big[\prod_{i=1}^m (1 &-V_i (X_{n,t_i})) \Big]
= \sum_{k_1 \leq \ldots \leq k_m} \P(N_{n,t_1}=k_1,\ldots, N_{n,t_m}=k_m)\\
& \quad \times \idotsint
\Big(1 - \prod_{i=1}^{m} (1-V_i (z_i) ) \Big)
\prod_{i=1}^{m} \big( p_{n, t_i - t_{i-1}}* \bar\rho_n^{(k_i-k_{i-1})*} \big) (z_i-z_{i-1})
\,\prod_{i=1}^m {\rm d}z_i,
\end{split}
\end{align}
where $\bar \rho_n^{k*}$ denotes the $k$-fold convolution of $\bar\rho_n$ with itself, and $(p_{n,t_i-t_{i-1}} *\bar\rho_n^{0*})(z):=p_{n,t_i-t_{i-1}}(z)$. We first rewrite
(\ref{eq:summandRep}) in a suitable form before applying Theorem~\ref{T:ri}.

On the RHS of (\ref{eq:summandRep}), for each $1\leq i\leq m$, we let $z_{i,1+k_i-k_{i-1}}:=z_i$ and rewrite
$$
\big( p_{n, t_i - t_{i-1}}* \bar\rho_n^{(k_i-k_{i-1})*} \big) (z_i-z_{i-1}) \, {\rm d}z_i
=  \idotsint p_{n, t_i - t_{i-1}}(z_{i,1}-z_{i-1})\!\!\! \prod_{j=1}^{k_i-k_{i-1}} \!\!\!\bar\rho_n(z_{i, j+1}-z_{i,j})\!\!\!\!\!\! \prod_{j=1}^{1+k_i-k_{i-1}}\!\!\!\!\!\!{\rm d}z_{i,j},
$$
as well as
$$
1-V_i(z_i) = \prod_{j=1}^{1+k_i-k_{i-1}} (1-V_{i,j}(z_{i,j})),
$$
where $V_{i,j}=0$ for all $1\leq j\leq k_i-k_{i-1}$ and $V_{i, 1+k_i-k_{i-1}}=V_i$. Integrating with respect to $\phi(z_0){\rm d}z_0$, the multiple integral in
(\ref{eq:summandRep}) is then in the same form as the LHS of (\ref{ri}), and therefore we can apply (\ref{ri}) to obtain

\begin{align*}
\idotsint
&\phi(z_0) \Big(1- \prod_{i=1}^{m} (1-V_i (z_i) ) \Big)
\prod_{i=1}^{m} \big( p_{n,t_i - t_{i-1}} * \bar\rho_n^{(k_i-k_{i-1})*} \big) (z_i-z_{i-1})
\,\prod_{i=0}^m {\rm d}z_i\\
&\geq
\idotsint
\phi_* (z_0) \Big(1- \prod_{i=1}^{m} (1-V_i^* (z_i) ) \Big)
\prod_{i=1}^{m} \big( p^*_{n, t_i - t_{i-1}} * (\bar\rho^*_n)^{(k_i-k_{i-1})*} \big) (z_i-z_{i-1})
\,\prod_{i=0}^m {\rm d}z_i.
\end{align*}
Since $C^*_{n,\cdot}$ has the same jump rate as $C_{n,\cdot}$ with jump kernel $\bar\rho_n^*$ instead of $\rho_n$, and $G^*_{n,\cdot}$ has transition
kernel $p^*_{n, t}$ (see e.g.~\cite[Sec.~3]{BM-H10}), summing the above inequality over $0\leq k_1\leq \cdots \leq k_m$ with weights $\P(N_{n,t_1}=k_1,\ldots, N_{n,t_m}=k_m)$
then gives \eqref{eq:approx}.
\qed

It was shown in the proof of~\cite[Theorem~4.3]{BM-H10} that $(X_{n,t_1},\ldots, X_{n,t_m})\Rightarrow(X_{t_1},\ldots, X_{t_m})$ and $(X^*_{n,t_1},\ldots, X^*_{n,t_m})\Rightarrow(X^*_{t_1},\ldots, X^*_{t_m})$ in distribution as $n\to\infty$. Using the Dominated Convergence Theorem, we can then easily extend Lemma~\ref{lem:gaussian} from $X_{n,\cdot}$ to $X$, which we state as follows.
\bp\label{prop:levy}
Let $X$ be a L{\'e}vy process with characteristic $(b, {\mathbb A}, \rho(x){\rm d}x)$. Let $\phi: \R^d\to [0,\infty)$, and let $m\in\N$. For $1\leq i\leq m$, let $V_i: \R^d\to [0,1]$ be continuous with compact support. Let $0 < t_1 < \ldots < t_m<\infty$. Then
\begin{align}\label{eq:levy}
\int_{\R^d} \phi(x) \Big(1- \E_x\Big[\prod_{i=1}^{m}(1-V_i(X_{t_i})) \Big] \Big) \, {\rm d}x \geq \int_{\R^d} \phi_*(x) \Big(1- \E_x\Big[\prod_{i=1}^{m}(1-V^*_i(X^*_{t_i})) \Big] \Big)\, {\rm d}x.
\end{align}
\ep
\bigskip

\noindent {\bf Proof of Theorem~\ref{T:trap}.} We may assume the conditions in Claim~\ref{clm:proofReduction}~(i)--(iii). Since $U$ is continuous with compact support and
$X$ is a.s.\ c{\`a}dl{\`a}g, for every $x\in\R^d$ and almost surely every realization of $X$ with $X_0=x$, we have
$$
\sum_{i=1}^{k}\frac{t}{k} U_{it/k} (-X_{it/k}) \asto{k} \int_0^t U_s(-X_s) \, {\rm d}s.
$$
By the same dominated convergence argument as in (\ref{domconv}), we have
\begin{align}\label{eq:expon1}
W^X_t(\phi, U_\cdot) = \lim_{k \to \infty} \int_{\R^d} \phi(x)
\Big(1 - \E_x\Big[ \exp \Big \{ -\sum_{i=1}^{k}\frac{t}{k} U_{it/k} (-X_{it/k}) \Big \} \Big] \Big) \, {\rm d}x.
\end{align}
By Lemma~\ref{L:cts*}, $U^*$ is also continuous with compact support. Therefore the same argument yields
\begin{align}\label{eq:expon2}
W^{X^*}_t(\phi_*, U^*_\cdot) = \lim_{k \to \infty} \int_{\R^d} \phi_*(x)
\Big(1 - \E_x\Big[ \exp \Big \{ -\sum_{i=1}^{k}\frac{t}{k} U^*_{it/k} (-X^*_{it/k}) \Big \} \Big] \Big) \, {\rm d}x.
\end{align}
Since for each $s=it/k$, we can write $e^{-U_s(-x)} = 1 - V_s(x)$ for a continuous $V_s: \R^d\to [0,1]$ with compact support, and note that
$1-V_s^*(x) = e^{-U_s^*(-x)}$, we can apply Proposition~\ref{prop:levy} combined with (\ref{eq:expon1})--(\ref{eq:expon2}) to obtain $W^X_t(\phi, U_\cdot)\geq W^{X^*}_t(\phi_*, U^*_\cdot)$.
\qed

\section{Proof of Theorem \ref{T:cap}}\label{S:cap}
We will deduce Theorem~\ref{T:cap} from Theorem~\ref{T:trap}.
\medskip

\noindent
{\bf Proof of Theorem~\ref{T:cap}.} Let $O$ be an open set, and recall that $T_O(X):=\inf\{s\geq 0: X_s\in O\}$. Note that applying Theorem~\ref{T:trap} with $\phi\equiv 1$ and $U_s(x)=\infty\cdot 1_O(-x)$ for all $s\geq 0$ gives
\be\label{capOcomp}
\int_{\R^d} \P_x(T_O(X)< t) \, {\rm d}x \geq \int_{\R^d} \P_x(T_{O^*}(X^*)< t) \, {\rm d}x \qquad \mbox{for all } t> 0.
\ee

We first consider $q>0$. By Definition~\ref{def:qCapa},
\be
\begin{aligned}
 C^q_X (O) &= q\int_{\R^d} \E_x\big[e^{- q T_O(X)}\big] \, {\rm d}x = q\int_{\R^d} \int_0^1 \P_x\big(e^{- q T_O(X)}> s\big) \, {\rm d}s \, {\rm d}x
 \\
&= q \int_0^1 \int_{\R^d} \P_x(T_O(X)< -q^{-1}\log s) \, {\rm d}x\,{\rm d}s \\
&\geq q \int_0^1 \int_{\R^d}  \P_x(T_{O^*}(X^*)< -q^{-1}\log s) \, {\rm d}x\,{\rm d}s \\
&= q\int_{\R^d} \int_0^1 \P_x\big(e^{- q T_{O^*}(X^*)}> s\big) \, {\rm d}s \, {\rm d}x = q\int_{\R^d} \E_x\big[e^{- q T_{O^*}(X^*)}\big] = C^q_{X^*} (O^*),
\end{aligned}
\ee
where in the inequality we applied (\ref{capOcomp}). For a general $A\in \Bi(\R^d)$, by (\ref{qChoquet}), we have
$$
C^q_X(A) = \inf\{C^q_X(O): A \subset O, \, O \mbox{ open}\}.
$$
Since for any open $O\supset A$, we have just proved that $C^q_X(O) \geq C^q_{X^*}(O^*)$, and $C^q_{X^*}(O^*)\geq C^q_{X^*}(A^*)$ by (\ref{capincl}),
we conclude that $C^q_X(A) \geq C^q_{X^*}(A^*)$ for all $A\in \Bi(\R^d)$. This proves Theorem~\ref{T:cap} for $q>0$.

Now consider the case $X$ is transient and $q=0$. By Definition~\ref{def:0Capa}, for any relatively compact $A\subset \R^d$,
\be
C^0_X(A) = \lim_{q\downarrow 0} C^q_X(A) \geq \lim_{q\downarrow 0} C^q_{X^*}(A^*) = C^0_{X^*}(A^*).
\ee
For general $A\in \Bi(\R^d)$, by definition, we have
$$
C^0_X(A) := \sup\{ C^0_X(K) : K\subset A,\, K\, {\rm relatively\ compact}\}.
$$
Let $A_n:= A\cap \{x \in\R^d: |x| \leq n\}$. Then $(A_n)_{n\in\N}$ are relatively compact, and
$$
C^0_X(A) \geq C^0_X(A_n) \geq C^0_{X^*}(A_n^*).
$$
Note that $(A_n^*)_{n\in\N}$ are finite open balls centered at the origin, and $A_n^* \uparrow A^*$ as $n\to\infty$, which implies that $C^0_{X^*}(A_n^*)\uparrow C^0_{X^*}(A^*)$ as $n\to\infty$ by (\ref{qChoquet}) and (\ref{capincl}). Therefore $C^0_X(A) \geq C^0_{X^*}(A^*)$ for all $A\in \Bi(\R^d)$, which proves Theorem~\ref{T:cap} for the case $X$ is transient and $q=0$.
\qed

\section{Some Open Questions}\label{sec:conclusion}
One of the open problems formulated at the end of~\cite{PS11} is the following. If $X$ is a standard Brownian motion with $X_0=0$, $f:[0,\infty)\to\R^d$ is measurable
(or even c\`adl\`ag or continuous), for which open sets $D$ of finite volume, is the expected volume of the Wiener sausage $\bigcup_{0\leq s\leq t}(D+X_s+f(s))$ minimized when we take $f\equiv 0$? For such $D$, then in light of the discussion after Remark~\ref{R:trap}, we will call the phenomenon where the optimal path is the
constant path, the {\em Pascal principle}. By the derivation leading to Corollary~\ref{C:sausage}, this question is equivalent to a trapping problem, where in Theorem~\ref{T:trap}, we take $\phi\equiv 1$, $U_s(x)=\infty\cdot 1_D(x-f(s))$, and ask whether $W^X_t(1, U_\cdot)$ is minimized at $f\equiv 0$. Note that because we
are not allowed to symmetrically rearrange $D$, standard rearrangement inequalities will not be applicable. Generalizing from Brownian motion, we may also ask if the above Pascal principle holds for any L\'evy process $X$ whose law is equally distributed with $X^*$.

In light of the analogy between the random walk exit problem in (\ref{survprob}) and the trapping problem in (\ref{exitmass}), we can ask whether the Pascal principle
holds for the survival probability of a Brownian motion killed upon exiting a finite domain. More precisely, let $X$ be a standard Brownian motion (or more generally a L\'evy process whose law is equally distributed with $X^*$), let $D$ be a closed set of finite volume with a sufficiently regular boundary, and assume that $X_0$ is distributed uniformly on $D$. Let $f:[0,\infty)\to\R^d$ be measurable (or even c\`adl\`ag or continuous), and let $T_{D^c}(X+f):=\inf\{s\geq 0: X_s+f(s)\in D^c\}$.
$$
\mbox{For which $D$ is} \quad \P(T_{D^c}(X+f)>t) \quad \mbox{maximized at } \ f\equiv 0 \mbox{ ?}
$$
When $D$ is symmetric and convex, $X$ is a standard Brownian motion, and $f$ is c\`adl\`ag, the answer is affirmative and it follows from Anderson's inequality~\cite{A55} for multi-variate normal distributions (Anderson's inequality is in fact valid for general symmetric unimodal distributions). If we do not impose any assumption on
the distribution of $X_0$, it is easily seen that the Pascal principle will fail in general. The uniform distribution on $D$ we propose is based on the analogy with the trapping problem. Another natural distribution for $X_0$ we may consider is the quasi-stationary distribution of $X$ on $D$, which equals the limit of $\P(X_t\in \cdot |T_{D^c}(X)>t)$ as $t\to\infty$.

\appendix

\section{Properties of Symmetric Rearrangements} \label{S:conv}
We collect here some basic properties of symmetric decreasing rearrangements that we use in the proof. Many facts here are standard to experts in rearrangement inequalities. However we include their proof for the sake of completeness, as well as for the convenience of the general reader. Below, Lemmas~\ref{L:approx} and \ref{lem:conv} are used to carry out approximations. Lemma~\ref{L:conv} considers the symmetric decreasing rearrangement of the convolution of two functions, while Lemma~\ref{L:cts*} considers the spatial symmetric decreasing rearrangement of a function which is continuous in space and time.

\blem\label{L:approx}
Let $\phi, \phi_n: \R^d\to [0,\infty]$, $n\in\N$, be such that $\phi_n(x)\uparrow \phi(x)$ as $n\to\infty$ for Lebesgue almost every $x\in\R^d$.
Then $\phi_n^*(x)\uparrow \phi^*(x)$ for every $x\in\R^d$.
\elem

\noindent
{\bf Proof.} Since $(\phi_n^*)_{n\in\N}$ and $\phi^*$ remain unchanged if $(\phi_n)_{n\in\N}$ and $\phi$ are modified on a set of Lebesgue measure $0$, we may assume without loss of generality that $\phi_n(x)\uparrow \phi(x)$ for all $x\in\R^d$. We now define for all $t>0$ the level sets
\be\label{PhinPhi}
\Phi_n(t) := \{z: \phi_n(z)>t \} \ \ \text{ and } \ \ \Phi(t):=\{z: \phi(z)>t \}.
\ee
Then by the assumption that $\phi_n \uparrow \phi$, we get that
\[
\Phi_n(t) \uparrow  \Phi(t)   \text{ as } n \to \infty,
\]
which implies that
\begin{align*}
\Phi^*_n(t) \uparrow \Phi^*(t) \text{ as } n \to \infty.
\end{align*}
(Note that if $|A|=\infty$, then we define $A^*:=\R^d$.) Therefore by the Monotone Convergence theorem,
\begin{align*}
\phi_n^*(x) = \int_{0}^{\infty} 1_{\Phi_n^*(t)}(x) \, {\rm d}t\
\big\uparrow\ \int_{0}^{\infty} 1_{\Phi^*(t)}(x) \, {\rm d}t \quad \text{ as } n \to \infty.
\end{align*}
Since $\phi^*(x) = \int_{0}^{\infty} 1_{\Phi^*(t)}(x) \, {\rm d}t$, we obtain
\[
\lim_{n \to \infty}\phi_n^*(x) = \phi^*(x).
\]
Note that this convergence hol{\rm d}s for every $x\in\R^d$.
\qed
\bigskip

\blem\label{lem:conv}
Let $\phi, \phi_n: \R^d\to [0,\infty)$, $n\in\N$, be uniformly bounded with uniformly bounded support, such that $\phi_n(x)\to \phi(x)$ as $n\to\infty$ for Lebesgue almost every $x\in\R^d$. Then $\phi_n^*(x)\to \phi^*(x)$ for Lebesgue almost every $x\in\R^d$.
\elem

\brm{\rm Lemma~\ref{lem:conv} is a correction of~\cite[Lemma 4.2]{BM-H10}, where $(\phi_n)_{n\in\N}$ were not assumed to have uniformly bounded support, and the conclusion can be seen to be false. Indeed, fix any $0\neq v\in\R^d$. Then $\phi_n(x) := 1_{\{|x-nv|< 1\}}$ converges pointwise to $\phi\equiv 0$, and yet $\phi_n^*(x)= 1_{\{|x|< 1\}}\not\to \phi^*\equiv 0$. }
\erm

\noindent
{\bf Proof.} As in the proof of Lemma~\ref{L:approx}, we may assume without loss of generality that $\phi_n(x)\to\phi(x)$ for every $x\in\R^d$. We will first show that $\liminf_{n \to \infty} \phi^*_n(x) = \phi^*(x)$ for Lebesgue a.e.\ $x$.

Since $(\phi_n)_{n\in\N}$ and $\phi$ are uniformly bounded with uniformly bounded support, by
the Dominated Convergence theorem, we have
\begin{align*}
\int_{\R^d} \phi_n(x) \, {\rm d}x \to \int_{\R^d} \phi(x)\, {\rm d}x \quad \text{ as } n \to \infty.
\end{align*}
Since for any nonnegative $f$, we have $\int f(x) \, {\rm d}x = \int f^*(x) \, {\rm d}x$, we obtain
\begin{align}\label{eq:limit}
\int_{\R^d} \phi_n^*(x) \, {\rm d}x \to \int_{\R^d} \phi^*(x)\, {\rm d}x \quad \text{ as } n \to \infty.
\end{align}
For $t>0$, let the level sets $(\Phi_n(t))_{n\in\N}$ and $\Phi(t)$ be defined as in (\ref{PhinPhi}). These level sets have finite volume by the assumption
that $(\phi_n)_{n\in\N}$ and $\phi$ have uniformly bounded support.

By the convergence of $\phi_n$ to $\phi$, we have
\begin{align*}
\Phi(t) \subset \liminf_{n \to \infty} \Phi_n(t) = \bigcup_{k=1}^{\infty} \bigcap_{n \geq k} \Phi_n(t).
\end{align*}
Using this, Lemma~\ref{L:approx} applied to the corresponding indicator functions yields the second equality in
\begin{align*}
\Phi^*(t) \subseteq  \Big(\bigcup_{k=1}^{\infty} \bigcap_{n \geq k} \Phi_n(t)\Big)^*
=\bigcup_{k=1}^{\infty} \Big( \bigcap_{n \geq k} \Phi_n(t)\Big)^* \subseteq \bigcup_{k=1}^{\infty}\bigcap_{n \geq k} \Phi_n^*(t) = \liminf_{n \to \infty} \Phi^*_n(t).
\end{align*}
We can now write
\begin{align*}
\phi^*(x) &= \int_{0}^{\infty} 1_{ \{ x \in \Phi^*(t) \} } \, {\rm d}t
\leq \int_{0}^{\infty} 1_{\liminf \Phi_n^*(t)}(x) \, {\rm d}t\\
&\leq \liminf_{n \to \infty} \int_{0}^{\infty} 1_{\Phi_n^*(t)}(x)\, {\rm d}t = \liminf_{n \to \infty} \phi_n^*(x),
\end{align*}
where in the second inequality we used Fatou's Lemma.
We thus showed $\displaystyle \phi^*(x) \leq \liminf_{n \to \infty} \phi_n^*(x)$ for all $x$. If we now integrate over all $x \in \R^d$ and use Fatou's lemma again, we obtain
\begin{align*}
\int_{\R^d} \phi^*(x) \, {\rm d}x \leq \int_{\R^d}\liminf_{n \to \infty} \phi_n^*(x) \, {\rm d} x
\leq \liminf_{n \to \infty} \int_{\R^d} \phi_n^*(x) \, {\rm d}x = \int_{\R^d} \phi^*(x) \, {\rm d}x,
\end{align*}
where the equality follows from \eqref{eq:limit}. Therefore, we deduce that
\begin{align}\label{eq:liminf}
\liminf_{n \to \infty} \phi_n^*(x) = \phi^*(x), \ \ \text{for Lebesgue a.e. } x.
\end{align}
We will now finish the proof by showing that $\limsup_{n \to \infty}\phi_n^*(x) \leq \phi^*(x)$ for Lebesgue a.e.\ $x$.
We define a new sequence of functions $f_n := \sup_{k \geq n} \phi_k(x)$.
By the assumptions on $(\phi_n)_{n\in\N}$, $(f_n)_{n\in\N}$ are also uniformly bounded with uniformly bounded support. Clearly,
$f_n \downarrow \phi$ as $n \to \infty$. Therefore, we may apply \eqref{eq:liminf} to $f_n$ instead of $\phi_n,$
and deduce that
\[
\liminf_{n \to \infty} f_n^*(x) = \phi^*(x) \ \ \text{ for Lebesgue a.e. } x.
\]
Note that $\liminf_{n\to\infty} f_n^*(x)=\limsup_{n \to \infty}f_n^*(x)$ because $f_n^*(x)$ is a nonincreasing sequence. Together with $\phi_n^*(x) \leq f_n^*(x)$ for all $x$, we obtain
\[
\limsup_{n \to \infty} \phi_n^*(x) \leq \phi^*(x) \ \ \text{ for Lebesgue a.e. } x,
\]
which concludes the proof.
\qed

\blem\label{L:conv} Suppose that $f,g: \R^d\to [0,\infty)$ and $f=f^*$, $g=g^*$. Then $f*g = (f*g)^*$.
\elem
{\bf Proof.} Since $f$ and $g$ are radially symmetric, so must be $f*g$. Since $f$ and $g$ are lower semi-continuous, for any $x_n\to x$, we have
$$
(f*g)(x) = \! \int\! f(x-y)g(y) \, {\rm d}y \leq \! \int\! \liminf_{n\to\infty}f(x_n-y) g(y) \, {\rm d}y \leq \liminf_{n\to\infty} \! \int\! f(x_n-y)g(y) \, {\rm d}y=\liminf_{n\to\infty} (f*g)(x_n).
$$
Therefore $f*g$ is also lower semi-continuous. It only remains to show that $f*g$ is radially nonincreasing. By writing
$$
(f*g)(x) = \int_{\R^d} \int_0^\infty 1_{\{f(x-y)>s\}} \,{\rm d}s  \int_0^\infty 1_{\{g(y)>t\}} \, {\rm d}t \, {\rm d}y
= \int_0^\infty \int_0^\infty (1_{F_s}*1_{G_t})(x) \,{\rm d}s\, {\rm d}t,
$$
where $F_s:=\{x: f(x)>s\}$ and $G_t:=\{x: g(x)>t\}$ are centered open balls, we only need to show that $(1_{F_s}*1_{G_t})(x)$ is radially nonincreasing.
This is equivalent to showing that $|F_s \cap (G_t+ \lambda x)|$ is nonincreasing in $\lambda\geq 0$ for any $x\neq 0$, which is clearly true.
\qed

\blem\label{L:cts*} Let $U_s(x) : [0,\infty)\times \R^d\to [0,\infty)$ be continuous with compact support. For each $s\geq 0$, let $U_s^*(\cdot)$ denote the symmetric decreasing
rearrangement of $U_s(\cdot)$. Then $U^*_\cdot(\cdot)$ is also continuous on $[0,\infty)\times\R^d$ with compact support.
\elem
{\bf Proof.} Clearly $U^*_\cdot(\cdot)$ has compact support. We first show that for each $s\geq 0$, $U^*_s$ is continuous. By definition of $U^*_s$, there exists
$f_s: [0,\infty)\to [0,\infty)$, which is non-increasing and right-continuous, such that $U^*_s(x)=f_s(|x|)$ for all $x\in\R^d$. If $U^*_s$ is discontinuous at
some $x_0\in\R^d$, then $f_s$ has a jump discontinuity at $|x_0|$, and we must have $|x_0|>0$. In particular, we must have
$$
0<|\{x\in \R^d: U_s(x)> f_s(|x_0|)\}| = |\{x\in \R^d: U_s(x)> f_s(|x_0|)+\eps\}|<\infty \qquad \mbox{for some } \eps>0.
$$
However the equality cannot hold because $U_s$ is continuous. Therefore $f_s$ must be continuous, and hence $U_s^*$ must be continuous as well.

Next we show that $f_s(r)$ is jointly continuous in $s\geq 0$ and $r\geq 0$. The continuity of $U_\cdot(\cdot)$ and Lemma~\ref{lem:conv} imply that for
each $s\geq 0$ and for Lebesgue a.e.\ $x\in\R^d$, $U^*_t(x)\to U^*_s(x)$ as $t\to s$. This in turn implies that for Lebesgue a.e.\ $r\geq 0$,
$f_t(r)\to f_s(r)$ as $t\to s$. Since $(f_s)_{s\geq 0}$ are all continuous, monotone, with uniformly bounded support, $f_t(\cdot)$ must converge uniformly
to $f_s(\cdot)$ as $t\to s$. This establishes the joint continuity of $f_s(r)$ in $s,r\geq 0$, and hence $U^*_s(x)$ must also be jointly continuous in
$s\geq 0$ and $x\in\R^d$.
\qed
\bigskip

\noindent
{\bf Acknowledgement} We thank Frank Aurzada for enlightening discussions on random walk exit problems and Yuval Peres for useful remarks. We are grateful to the referee for providing many helpful comments and references. R.~Sun is supported by grant R-146-000-119-133 from the National University of Singapore and
A.~Drewitz has been supported by an ETH Fellowship.

\end{document}